\documentclass[12pt,a4paper]{article}
\usepackage{amsmath}
\usepackage{graphicx}
\usepackage{amssymb}
\usepackage{hyperref}
\begin{document}
\title{Theorem of Wantzel}
\author{J. Mainik}
\date{June 20, 2026}
\maketitle
\begin{abstract}
In 1796, Gauss succeeded in solving the problem of constructing the regular
17-gon with compass and straightedge. Later he proved that, using a compass and straightedge,  it is possible
to construct the regular polygons with
$n=2^m n_1\cdots n_l$ sides if $n_1,\cdots, n_l$ are different prime
numbers of the form $n_k=2^{2^{\nu_k}}+1$.
Gauss also knew that only these regular polygons can be constructed but did
not prove it.\linebreak
P. Wantzel completed the result of Gauss and proved it in 1837.
The present paper provides a new proof for Wantzel's theorem.
\end{abstract}
 

\section{Introduction}

The constructing of regular polygons with a compass and straightedge
is a well-known task that has been engaging mathematicians for over 2000 years.
In 1796, C. F. Gauss found first the solution for the $17$-gon,
\cite{Gauss1}-\cite{Gauss2},
and proved later that it is possible to construct
 with a compass and straightedge
the regular polygons with $n=2^m n_1\cdots n_l$ sides, where $n_1,\cdots, n_l$ are different prime
numbers of the form $n_k=2^{2^{\nu_k}}+1$.

Gauss also knew that only these regular polygons can be constructed.
He sought for a proof but has not finished 
this work. It may be that he was too busy with other tasks or, due to the zeitgeist
of the era before the successful development of the Galois theory,
the interest in the proof of impossibility was conceptually not 
yet well-established.

P. Wantzel completed the result of Gauss and proved the impossibility in 1837,
\cite{Wantzel}.
Because of his early death 
Wantzel's proof  
and Wantzel himself were unfortunately forgotten for almost a century.
In the brilliant publication by F. Klein
\cite{Klein}
in 1895 a new proof was presented, and there was no reference to
Wantzel and his publication.
A repetition of Wantzel's proof and a detailed description of the 
related matters can be found in the paper 
\cite{Lützen}
from 2009.

Here we present a new proof of Wantzel's theorem, and a part of this
proof will be done using two different approaches.
A justification for this work and, in particular, for the two approaches 
in a part of the proof, so author believes, can be given using the remark of Gauss
that "{\it it is not out of place to discover 
the same truths by different methods}",
\cite{Gauss2}.
The author hopes that the new proof of Wantzel's theorem will find interest
and understanding.

The numbers of the form $n=2^{2^{\nu}}+1$ are named Fermat numbers.
Fermat thought that they are prime numbers, but this is not correct.
The numbers
$2^{2^0}+1=3, 2^{2^1}+1=5, 2^{2^2}+1=17, 2^{2^3}+1=257$ and $2^{2^4}+1=65537$
are indeed primes, but the numbers $2^{2^\nu}+1 $, $5\le \nu\le 32$,
are not.
It is currently unknown whether other Fermat primes exist.
This problem is not trivial:
the smallest Fermat number for which it is not known whether it is prime or not is $2^{2^{33}}+1$, and this is a number with $2\,585\,827\,973$ digits.

It is clear that the regular $n$-gon with $n=2^m n_1\cdots n_l$ sides,
where $n_1,\cdots ,n_l$ are different primes, 
can be constructed exactly when every single $n_j$-gon can be constructed.
Therefore the proof of Wantzel's theorem consists of the following two parts:

1) Show that the regular $n$-gon cannot be constructed
if the prime $n$ is not a Fermat number.

2) Show that for a Fermat prime $n$ the regular $n^2$-gon
cannot be constructed. It is sufficient to analyze only this case.
If we could construct the $n^k$-gon, then we could also construct the $n^2$-gon.

\section{Denominations and remarks}

For $z=e^{2\pi i/n}$ we have 
\begin{eqnarray}\label{Gl1}
z^n-1=0,
\end{eqnarray}
and the values $e^{m\cdot 2\pi i/n}$, $0\le m\le n-1$, are the solutions
of this equation. We denote $v=e^{2\pi i/n}$ and present these values
as $v^m$, $0\le m\le n-1$. 
These points are located on the circle with radius $1$ and represent the vertices 
of the regular $n$-gon.
In the following we will refer to these values as \emph{elements}.

It holds obviously
\begin{eqnarray}\label{SummexHochk}
z^n-1=(z-1)(z^{n-1}+z^{n-2}+\cdots +z+1).\nonumber
\end{eqnarray}
We have the radius of the circle equal to $1$ and the point $v^0=1$.
The values $v^m$, $1\le m\le n-1$, are the solutions of the
cyclotomic equation:
\begin{eqnarray}\label{Kreisteilung}
z^{n-1}+z^{n-2}+\cdots +z+1=0.
\end{eqnarray}
To construct the regular $n$-gon, we have to determine the solutions
of this equation and they, in addition, must be constructed with a compass 
and straightedge. The polynomial 
\begin{equation}\label{KreisPolyn}
P(z)=z^{n-1}+z^{n-2}+\cdots +z+1
\end{equation}
is called cyclotomic polynomial, and it is easy to see from 
\eqref{KreisPolyn}
that the sum of the roots of this polynomial is equal to $-1$: 
$v^1+v^2+\cdots +v^{n-1}=-1$.

It is a natural question: what numbers (i.e. lengths of segments) are we 
dealing with when constructing
with a compass and straightedge. In any case 
we have the rational numbers because it is possible
to construct the segment of length $m$ and divide it into $n$ parts
of the same length $m/n$. But they are not enough.
If we consider the points at which a line and a circle or two circles intersect,
we get a quadratic equation, and we need the square root $\sqrt D$ of the
corresponding determinant. 
So we come to extensions of the rational numbers involving
the adjunction of square roots.

At the first 
extension of the rational numbers $\mathbb{Q}=\mathbb{Q}_0$ with
the square
$\sqrt{D_1}$, where  $D_1\in \mathbb{Q}$  and $\sqrt{D_1}\not \in \mathbb{Q}$,
we get the linear combinations $a+b\sqrt{D_1}$ with $a\in \mathbb{Q}$ and $b\in \mathbb{Q}$
and these numbers build a number field $\mathbb{Q}_1=\mathbb{Q}(\sqrt{D_1})$.
It is clear how the arithmetic operations for these numbers 
can be performed. 

These extensions can be repeated a finite number of times.
At the $m$-th extension we thus get a field $\mathbb{Q}_m$ of numbers
$a+b\sqrt{D_m}$ for which $a$, $b$ and $D_m$ belong to the previous
extension $\mathbb{Q}_{m-1}$ and $\sqrt{D_m}$ does not belong to $\mathbb{Q}_{m-1}$.
It is clear again how to perform the arithmetic operations.
The equality $a_1+b_1\sqrt{D_m} = a_2+b_2 \sqrt{D_m}$
holds here exactly when $a_1=a_2$ and $b_1=b_2$.

If we consider the intersections of straight lines 
and the distances between points, no new demands will arise. 

The segment of length  $a+b\sqrt D$ itself
can be constructed for given segments of length  $a$, $b$ and $D$ 
since we have the segment of length $1$.
\newline
\newline
In our proof we will use the following known results:

1. Every prime number has primitive roots,
\cite{Gauss1}.

An integer $g$ is a primitive root for the modulo $n$ if the remainder
$\mathrm{rest}(g^k;n)$, $1\le k\le n-1$, of the division $g^k $ by $n$
results in all values $1,2,\cdots,n-1$.

2. If $n$ is a prime number, then the corresponding cyclotomic polynomial
$P(z)$ is irreducible, i.e. the polynomial $P(z)$ cannot be 
factored into two polynomials with rational coefficients.

The proof of this result was obtained as follows.
First, Gauss proved the following lemma: If a polynomial 
\[
F(z)=z^k +c_1z^{k-1}+c_2z^{k-2}+\cdots+c_k
\]
with integer coefficients can be factored into the product of two
polynomials with rational coefficients
\[
F(z)=(z^m +\alpha_1z^{m-1}+\cdots \alpha_m)\cdot(z^{k-m}+\beta_1z^{k-m-1}+\cdots \beta_{k-m})
\]
then this polynomial can be factored into the product of two polynomials
with integer coefficients.
Later, Eisenstein proved that the cyclotomic polynomial cannot be 
factored into the product of two polynomials with integer coefficients,
\cite{Eisenstein}.
A detailed proof can be found in 
\cite{Klein}.

3. Let $F(x)=x^k +c_1x^{k-1}+c_2x^{k-2}+\cdots+c_k$ be a polynomial
of odd degree with integer coefficients.
If this polynomial does not have an integer root, then not all real roots 
of this polynomial can be constructed with a compass and straightedge.

One arrives at the proof as follows.

a) It is easy to see that the polynomial $F(x)$ has no 
rational
roots that are not integers. The roots of $F(x)$ 
must therefore be determined using a finite number of quadratic extensions 
of the rational numbers.

b) If, in the $m$-th extension of the rational numbers, the value\linebreak
$x_1=a+b\sqrt{ D_m}$ is a root of the polynomial $F(x)$ with
coefficients from the $(m-1)$-th extension, then the value $x_2=a-b\sqrt{D_m}$ 
is also a root of this polynomial. Indeed, 
for the polynomial with coefficients from the $(m-1)$-th extension 
we have
$F(x_1)=F(a+b\sqrt {D_m})=a_1+b_1\sqrt{D_m}=0$, where $a_1$ and $b_1$ 
belong to the $(m-1)$-th extension and therefore $a_1=0$ and $b_1=0$.
For  $x_2$ we get here $F(x_2)=F(a-b\sqrt{D_m})=a_1-b_1\sqrt{D_m}$,
and $F(a-b\sqrt{D_m})=0$ because of $a_1=0$ and $b_1=0$.

c) If we divide $F(x)$ by $(x-x_1)(x-x_2)=x^2-2ax+(a^2-b^2D_m)$, we
get a polynomial $F_1(x)$ with coefficients from the $(m-1)$-th extension.
Compared to $F(x)$, the polynomial $F_1(x)$ has lost two roots.
One can continue in this way and thus see that for a multiple root $a+b\sqrt{D_m}$ of $F(x)$  the value $a-b\sqrt{D_m}$ ist also a root
exactly as often.

It follows from a), b) and c) that 
$F(x)$ 
has an even number of real roots $x_1,\cdots, x_{2k_0}$ that can be constructed
using a compass and straightedge and that division $F(x)$ by
$(x-x_1)\cdots (x-x_{2k_0})$ yields a polynomial $F_0(x)$ that has 
no constructible real roots. However, $F_0(x)$ is a polynomial of odd degree 
and has at least one additional real root, which is also a root
of $F(x)$.
\newline

Finally, a few remarks on primitive roots. 
According to Fermat’s little theorem, $g^{n-1}\equiv 1 \pmod{n}$ if
$g \not\equiv 0\pmod{n}$.
Therefore we can in the definition of the primitive roots,
instead of $g^1,g^2,\cdots, g^{n-1}$, use the values
$g^0=1,g^1,\cdots, g^{n-2}$.
It is clear that $g$ is a primitive root if and only if $n-1$
is the smallest number for which $g^{n-1} \equiv 1 \pmod{n}$.

The following simple property will be useful:
If $g$ is a primitive root for the prime number $n$, then 
it must hold
$g^{\frac{n-1}{2}}\equiv -1\pmod{n}$. 
It follows then $g^{\frac{n-1}{2}+k}\equiv -g^k\pmod{n}$ for $k\ge 1$,
which means that
\begin{eqnarray}\label {inverseNr}
\mathrm{rest}(g^{\frac{n-1}{2}+k-1};n)=n-\mathrm{rest}(g^{k-1};n), 1\le k\le\frac{n-1}{2}.
\end{eqnarray}

\section{Primes that are not Fermat numbers}

{\bf Proposition 1.} If the prime number $n$ is not a Fermat number,
then the regular $n$-gon cannot be constructed using a compass and straightedge.
\newline
\newline
{\bf Proof.} Suppose the prime number $n > 3$ is not a Fermat number. Then  
the even number $n-1$ has an odd divisor $k$, $k > 1$, so that
$(n-1) = n_0 \cdot k$. Using the primitive root $g$ for $n$, 
we form $n-1$ distinct numbers $r_k = \mathrm{rest}(g^{k-1}; n)$, 
and split the ordered set $V$ of these numbers 
$V = \{r_1, r_2, \cdots, r_{n-1}\}$ into $k$ disjunct and ordered 
parts $V_j$, $1 \le j \le k$,
\[
V_j=\{r_j,r_{j+k},r_{j+2k},\dots,r_{j+n-1-k}\}
\]
Thus, $V_j$ contains the $j$-th number and every $k$-th number thereafter.
For the sets $V_j$ defined in this way, 
the following cyclic property clearly holds. If we multiply the numbers
of the sets by the factor $g$ and calculate modulo $n$,
we obtain for $1\le j<k$ from $V_j$ the set $V_{j+1}$ and from $V_k$ 
the set $V_1$.

Each set $V_j$ consists of $n_0 = (n-1)/k$ numbers, and it follows from 
\eqref{inverseNr}
that the set of the $n_0$ values 
$v^{r_j},v^{r_{j+k}},v^{r_{j+2k}},\cdots,v^{r_{j+(n_0-1)\cdot k}}$ 
with 
degrees from $V_j$ consists of pairs of inverse elements.
The value $v^{r_j+(m-1)k}$ 
for the $m$-th number from the first half in $V_j$ 
(i.e. for $1\le m\le n_0/2$)
is inverse to the value $v^{r_{L/2+j+(m-1)k}}$ for the $(n_0/2+m)$-th number 
from the second half
in $V_j$.

We define the values $x_j$, $1\le j\le k$, as follows:
\begin{equation}\label{x_j}
x_j=\sum_{m\in V_j} v^m=\sum_{0\le m<n_0} e^{r_{j+km}\cdot 2\pi i/n}.
\end{equation}
Here and in what follows, we use the agreed notation $v=e^{2\pi i/n}$.
Due to the property of the values $v^m$ of degree $m$ in $V_j$ it immediately follows
that the values $x_j$ are real numbers,
\begin{equation}\label {reelx_j}
x_j=\sum_{0\le m< n_0/2 }2\cos (r_{j+km}\cdot 2\pi/n).
\end{equation}
In
\eqref {reelx_j}
only the first half of the numbers in $V_j$ is used.

We now consider the polynomial $F(x)=(x-x_1)(x-x_1)\cdots (x-x_k)$.
This is a polynomial of odd degree $k$ with real roots.
$F(x)$ has the form
\begin{equation}
F(x)=x^k+c_1x^{k-1}+c_2x^{k-2}+\cdots +c_{k-1}x+c_k,
\end{equation}
and we want to show that the coefficients $c_p$, $1\le p\le k$,
are integers. 
The coefficient $c_1$ is $1$ here:
\begin{equation}
c_1=-(x_1+\cdots +x_k)=-\sum_{1\le m\le n-1} v^m =1,
\end{equation}
and hence only the coefficients $c_p$, $2\le p\le k$, need to be considered.

For integers $k$ and $p$ such that $2 \le p \le k$ let $M(k,p)$ denote
the set of all distinct sequences $\{j_1, j_2, \cdots, j_p\}$ of integers, 
formed in such a way that\linebreak
 $1 \le j_1 < j_2 < \cdots < j_p \le k$.
For $c_p$
the following representation clearly holds:
\begin{equation}\label{c_m}
c_p = (-1)^p \sum_{M(k,p)} x_{j_1}x_{j_2}\cdots x_{j_p}
\end{equation}
and due to
\eqref{c_m} and \eqref{x_j}
we obtain the representation of $c_p$ in terms of the elements 
$v^1,\cdots, v^{n-1}$
\begin{equation}\label{cp_1}
c_p=\mu_{p,0}+\mu_{p,1}v^{r_1}+\mu_{p,2}v^{r_2}+\cdots \mu_{p,n-1}v^{r_{n-1}}
\end{equation}
with integer coefficients $\mu_{p,0}, \mu_{p,1}, \dots, \mu_{p,n-1}$.
The constant $\mu_{p,0}$ can be obtained here using products 
of inverse elements, due to $v^l\cdot v^{n-l}=1$.

We now multiply  the numbers in each set $V_j$, $1 \le j \le k$,
by the factor $g$ (calculating modulo $n$). 
With this change of the degrees in
\eqref{c_m}
for each factor in the product $x_{j_1}\cdots x_{j_p}$ the degrees $v^m$ 
of the corresponding summands in \nolinebreak 
\eqref{x_j}
are multiplied by $g$. 
If the degrees in the summands of the individual factors
$x_{j_1}, x_{j_2}, \dots, x_{j_p}$ are multiplied by the factor $g$, 
it is clear that in the result 
for the product $x_{j_1}x_{j_2}\cdots x_{j_p}$ 
the degrees $v^m$ of the corresponding summands are also multiplied by the factor $g$. 
Thus it follows that in the representation for $c_p$ in 
\eqref{cp_1}
all degrees of $v$ must be multiplied by the factor $g$.
This means that in 
\eqref{cp_1}
for $1\le j\le n-2$, instead of $\mu_{k,j}v^{r_j}$, we must get $\mu_{k,j}v^{r_{j+1}}$ 
and, instead of $\mu_{k,n-1}v^{r_{n-1}}$, we must get $\mu_{k,n-1}v^{r_{1}}$.
However, when multiplying the degrees by the factor $g$, 
we are going from $x_1, x_2, \dots, x_k$ to
$x_2, x_3, \dots, x_k, x_1$, and since $c_p$ does not depend on the order
of the roots $x_j$, $1 \le j \le k$, we must obtain the same value $c_p$.
Thus, we have
\begin{equation}\label{cp_2}
c_p=\mu_{p,0}+\mu_{p,1}v^{r_2}+\mu_{p,2}v^{r_3}+\cdots \mu_{p,n-2}v^{r_{n-1}}+\mu_{p,n-1}v^{r_1}.
\end{equation} 

Here it is appropriate to explain how the different representations
of $c_p$ should be understood. In what follows, we will need to compare analogous
representations on several occasions.
When calculating the coefficient $c_p$ in 
\eqref{cp_1},
the values
$x_j=v^{r_j}+v^{r_j+k}+\cdots+v^{r_j+(n_0-1)\cdot k}$ were used, and when calculating
$x_{j_1}x_{j_2}\cdots x_{j_p}$, one had to compute products of the form 
$v^{m_1}\cdot v^{m_2}$.
Since $v^{m_1}\cdot v^{m_2}=v^{m_1+m_2}$ one simply had to add the corresponding degrees. The coefficient $\mu_{p,k}$ in the representation
\eqref{cp_1}
is an integer that indicates how many times the element $v^{r_k}$ 
was finally obtained.

When calculating the coefficient 
$c_p$ in 
\eqref{cp_2},
we use the values\linebreak
$x_j=v^{r_{j+k}}+\cdots+v^{r_j+(n_0-1)\cdot k}+v^{r_j}$, which are obtained by 
multiplying the degrees by the factor $g$.
The order of the summands in $x_j$ has changed, but we still obtain 
the same products $x_{j_1}x_{j_2}\cdots x_{j_p}$ -- possibly,
with a different order of the summands.
No transformation is applied to the resulting sum of elements $v^{r_k}$, and
the quantity of summands $v^{r_k}$ (i.e. the coefficient $\mu_{p,k}$) cannot change in the process.

The two representations thus simply show that $\mu_{p,1}=\mu_{p,2}$,\linebreak 
$\mu_{p,2}=\mu_{p,3},\cdots,\mu_{p,n-2}=\mu_{p,n-1}$  and $\mu_{p,n-1}=\mu_{p,1}$,
because in the first calculation we have for $v^{r_{k+1}}$ the coefficient
$\mu_{p,k+1}$ 
and in the second calculation the coefficient $\mu_{p,k}$, 
and they must be the same.

The constant $\mu_{k,0}$ does not change because the inverse elements $v^l$
and $v^{n-l}$ become inverse elements when the degrees are multiplied by $g$.
Because of
\eqref{cp_1} and \eqref{cp_2}
it follows that in fact all coefficients $\mu_{p,m}$, $1\le m\le n-1$, are equal, 
and, since $v^1+v^2+\cdots +v^{n-1}=-1$, the following 
holds true:
\begin{equation}\label{cp}
c_p=\mu_{p,0}+\mu_{p,1}\cdot (v^1+v^2+\cdots +v^{n-1})=\mu_{p,0}-\mu_{p,1}.
\end{equation}
Thus, when calculating $c_p$, we obtain an even 
coverage of the elements $v^1,\cdots, v^{n-1}$, and the value $\mu_{p,1}$ 
is the height of this even coverage.
We have proved that $F(x)$ is a polynomial with integer coefficients.

Let the root $x_j$ of the polynomial $F(x)$ be an integer $C$.
Then consider the following polynomial $P_1(z)$:
\begin{equation}\label{P_1}
P_1(z)=\sum_{m\in V_j} z^m -C.
\end{equation}
$P_1(z)$ is a polynomial with integer coefficients, its degree is no greater than\linebreak
 $n-1$, and it is not equal to the cyclotomic polynomial $P(z)$.
For the polynomials $P(z)$ and $P_1(z)$
we have: $P(v^1)=0$ and $P_1(v^1)=0$.

We divide $P(z)$ by $P_1(z)$ and obtain
\begin{equation}\label{Divid}
P(z)=P_1(z)\cdot R_1(z)+P_2(z),
\end{equation}
where $P_2(z)$ is a polynomial with rational coefficients and 
of lower degree than $P_1(z)$.
Since $P(z)$ is irreducible, $P_2(z)$ is not $0$. 
For the value $v^1$ we obtain from
\eqref{Divid}
that $P_2(v^1)=0$, and one can divide $P(z)$ by
$P_2(z)$ and, by analogy, obtain a polynomial $P_3(z)$ of lower degree
with rational coefficients, for which $P_3(v^1)=0$ holds. 

In this way, we could determine polynomials of lower degree an infinite number of times, 
and this apparent contradiction
shows that the polynomial $F(x)$ has no integer roots. 
It immediately follows that the regular $n$-gon cannot be
constructed using a compass and a straightedge. Indeed, if we had the regular
$n$-gon,
we could construct all values $\cos(2\pi  m/n)$, $1\le m\le n-1$, 
and, using
\eqref{reelx_j},
also all values $x_j$. This completes the proof.
$\blacksquare$ 

{\bf Remark.} Apart from $c_1=1$, the coefficient $c_2$ can also be 
easily calculated. It applies here $c_2=-n_0\cdot (k-1)/2$. But, unfortunately, the other coefficients cannot be calculated so easily.
\section{Square of a Fermat Prime. Proof 1.}
We present two proofs for this case.
The first proof uses well-known results about the irreducibility of polynomials.

First, a few simple additional remarks and denominations.
For the $n^2$-gon,\linebreak 
instead of
\eqref{Kreisteilung},
we get the equation
\begin{equation}\label{nn-Division}
z^{n^2-1}+z^{n^2-2}+\cdots+z+1=0.
\end{equation}
We denote the polynomial in this equation by $P_{n^2-1}(z)$ and set 
$w=e^{2\pi i/n^2}$ here.
The roots of $P_{n^2-1}(z)$
are $w^m$, $1\le m\le n^2-1$, and the values $w^{mn}=v^m$, $1\le m\le n-1$, 
are the known roots of the cyclotomic polynomial $P(z)$ for the $n$-gon.

The polynomial $Q(z)$, with these $n(n-1)$ new roots, is given by
\begin{equation}\label{Q_nn}
Q(z)=P_{n^2-1}(z)/P(z)=z^{n(n-1)}+z^{n(n-2)}+\cdots+z+1
\end{equation}
and we see from 
\eqref{Q_nn}
that the sum of the $n(n-1)$ new roots $w^m$ with $m\neq kn$ is equal to $0$. 
We will denote $Q(z)$ the adapted cyclotomic polynomial.
In the proof of Wantzel’s theorem for primes that are not Fermat numbers
we used the property that the 
corresponding 
cyclotomic polynomial is irreducible. 
This 
also holds for the adapted cyclotomic polynomial.
The proof can be found in
\cite{Klein}.
In the following Proof 1 for the $n^2$-gon we will use this property.
\newline
\newline
{\bf Proposition 2.} If $n$ is a Fermat prime,
then the regular $n^2$-gon cannot be constructed with a compass and straightedge.
\newline
\newline
{\bf Proof 1.}
Let $g$ be a primitive root for the Fermat prime $n$. We denote by $U$ and $W$
the sets of the numbers
between $1$ and $n^2-1$ that are divisible by $n$
and non divisible by $n$, respectively:
$U=\{ n,2n,\dots,n(n-1)\}$,
$W =\{1,2,\dots,n^2-1 \}\backslash U$.

The numbers from $U$ are obtain (calculated modulo $n^2$) using 
$ng^{m-1}$
\[
\mathrm{rest}(ng^{m-1};n^2)=n\cdot \mathrm{rest}(g^{m-1};n), 1\le m\le n-1,
\]
so that $U=\{n\cdot r_1,n\cdot r_2,\dots ,n\cdot r_{n-1}\}$, where
as above $r_m=\mathrm{rest}(g^{m-1};n)$.

The numbers $g^{m-1}$, $1 \le m \le n(n-1)$, (calculated modulo $n^2$)
belong to \nolinebreak$W$. 
Indeed,
$g^{m-1} \not \in W$ means that $g^{m-1} \equiv 0 \pmod{n}$, and this cannot hold for the
primitive root $g$ for $n$.
However, it is not yet clear which numbers in $W$ can in fact be determined using
$g^{m-1}$. Let’s analyze this.

We have $g^{n-1} \equiv 1 \pmod{n}$, so that $g^{n-1} = jn + 1$ and we obtain
\[
g^{(n-1) n} = (jn + 1)^n = (jn)^n + n \cdot (jn)^{(n-1)} + \frac{n(n-1)}{2}
(jn)^{n-2}+\cdots+n\cdot jn +1
\]
and therefore $g^{n(n-1)}\equiv 1\pmod{n^2}$.

We have here $g^{\frac{n-1}{2}} \equiv -1 \pmod{n}$, so that
$g^{\frac{n-1}{2}} = jn - 1$. Then, by analogy with the previous result,
we conclude that $g^{\frac{n-1}{2}\cdot n}\equiv -1\pmod{n^2}$. 

The smallest period $p$ for which we have $g^p \equiv 1 \pmod{n^2}$ must be a divisor
of $n(n-1)$. 
This cannot be $p=\frac{(n-1)}{2^l}\cdot n$ with $l\ge 1$, 
since then
$g^{\frac{n-1}{2}\cdot n}\equiv 1\pmod{n^2}$. 
Nor can it be
$p=\frac{(n-1)}{2^l}$ with $l\ge 1$, since then 
$g^{\frac{n-1}{2}}\equiv 1\pmod{n^2}$ and thus also 
$g^{\frac{n-1}{2}\cdot n}\equiv 1\pmod{n^2}$.
Therefore, only two cases are possible: case a) with $p=(n-1)n$ and case b) with $p=n-1$.

In case a) all numbers in $W$ are determined by $g^{m-1}$, $1 \le m \le n(n-1)$,
and with $s_m = \mathrm{rest}(g^{m-1}; n^2)$, $1 \le m \le n(n-1)$, we have
an ordered set of all numbers in $W$. 
In this case, the proof can be conducted in exactly the same way as for prime numbers 
that are not Fermat numbers. Here, $W$ consists of disjunct ordered subsets
$W_j, 1\le j\le n$,
\[
W_j=\{s_j,s_{j+n},s_{j+2n},\dots,s_{j+n(n-2)}\}
\]
The values $x_j$, $1 \le j \le n$, are given by
\begin{equation}\label{xj_nn}
x_j = \sum_{m \in W_j} w^m 
\end{equation}
and they are again real numbers.
Indeed, we can show here, analogous to the case of a prime number, that
the numbers in the first and second half of $W$ are degrees of inverse elements.
Then it's the same for every $W_j$, and the corresponding set of elements
with degrees from $W_j$ consists of inverse pairs.

Here $F(x)=(x-x_1)(x-x_1)\cdots (x-x_n)$ 
is the required polynomial of odd degree $n$ with real roots.
$F(x)$ has the form
\begin{eqnarray}\label{F(x)_nn}
F(x)=x^n+c_1x^{n-1}+c_2x^{n-2}+\cdots +c_{n-1}x+c_n,
\end{eqnarray}
and for the coefficient $c_p$, $1\le p\le n$, we have the representation 
\eqref{c_m},
where $k$ must be set equal to $n$.
However, with the product $x_{j_1}x_{j_2}\cdots x_{j_p}$ 
one can obtain here the values $w^{n\cdot r_k}$
with $n\cdot r_k\in U$ and also the values $w^{s_k}$ with $s_k\in W$.
Instead of
\eqref{cp_1},
the following holds here:
\begin{eqnarray}\label{cp_nn_1}
c_p&=&\mu_{p,0}+\mu_{p,1}w^{n\cdot r_1}+\mu_{p,2}w^{n\cdot r_2}+\cdots \mu_{p,n-1}w^{n\cdot r_{n-1}} \nonumber \\
&+ &\delta_{p,1}w^{s_1}+\delta_{p,2}w^{s_2}+\cdots \delta_{p,n(n-1)}w^{s_{n(n-1)}}.
\end{eqnarray}

If we multiply the numbers in $W$ by the factor $g$ (computed modulo $n^2$),
we obtain the following result:
\begin{eqnarray}\label{cp_nn_2}
c_p&=&\mu_{p,0}+\mu_{p,1}w^{n\cdot r_2}+\mu_{p,2}w^{n\cdot r_3}+\cdots \mu_{p,n-2}w^{n\cdot r_{n-1}}+\mu_{p,n-1}w^{n\cdot r_1} \nonumber \\
&+&\delta_{p,1}w^{s_2}+\delta_{p,2}w^{s_3}+\cdots \delta_{p,n(n-1)-1}w^{s_{n(n-1)}}
+\delta_{p,n(n-1)}w^{s_1}.
\end{eqnarray}
It follows from
\eqref{cp_nn_1} and \eqref{cp_nn_2}
that all coefficients $\mu_{p,m}$, $1 \le m \le n-1$, 
and also all coefficients $\delta_{p,m}$, $1 \le m \le n(n-1)$, are equal, 
and the following representation holds:
\begin{equation}
c_p=\mu_{p,0}+\mu_{p,1}
\sum_{m\in U} w^m+ \delta_{p,1}\sum_{m\in W} w^m =\mu_{p,0}-\mu_{p,1},
\end{equation}
since $\sum_{m\in U}w^m=-1$ and $\sum_{m\in W}w^m=0$. 
If we calculate $c_p$ here, we obtain a constant and even coverages 
of the corresponding sets of elements,
and we need only the constant and the coverage height for 
$w^{n},\cdots,w^{(n-1)n}$.

After that, in case a), we can complete the proof in exactly the same way as in the case
of primes that are not Fermat numbers. We only have to use the irreducibility 
of the adapted cyclotomic polynomial $Q(z)$.\newline

Now let us consider case b). Here, using $g^{m-1}$, $m=1,\dots$ 
(computed modulo $n^2$), we can determine only $n-1$ numbers and not all
$n(n-1)$ numbers in $W$. The approach used above is not directly applicable here. 

In 
\cite{Mainik}
the author considered so-called invariant sets, and they were very helpful 
for a new method of constructing the regular polygons for Fermat primes. 
Here we want to use an invariant set in $W$ to build the appropriate
parts $W_j$.
An invariant set with the starting number $s\in  W$ for modulo $n^2$
consists of all distinct numbers  $\mathrm{rest}(s\cdot2^m;n^2)$, $m=0,1,\dots$
It is clear that the formation of the invariant set stops
when one returns to the starting number.
The invariant sets have the following nice circular property: 
After the last number in the set we return to the starting number 
and go through the same values again in the same order. 
We want to use the invariant set $H_1$ with the starting number $1$,
and we must examine it more closely.

The values $2^m$, $1\le m\le 2^{{\nu}}$, are between $2$ and $n-1$, and for $m=2^{\nu}$
we have $2^m=2^{2^{\nu}}=n-1$. Therefore, for $m=2^{\nu}\cdot s$ with $s\ge 1$ we have
\[
2^{2^{\nu}\cdot s}=(n-1)^s=n^s-s\cdot n^{s-1}+\frac{s(s-1)}{2}n^{s-2}+
\cdots+(-1)^{s-1} s n+(-1)^s
\]
and thus $2^{2^{\nu}\cdot s}\equiv (-1)^{s-1} sn+(-1)^s\pmod{n^2}$.
Then it is immediately clear that $s=2n$ is the smallest number for which 
$2^{2^{\nu}\cdot s}\equiv 1\pmod{n^2}$ holds.

For the values $2^{m}$ (calculated modulo $n^2$) we return to the starting 
number \nolinebreak$1$ when $m=n\cdot 2^{\nu+1}$,
and we do not return to $1$ if $m=2^{\mu}$ with $1\le \mu\le \nu+1$.

We now want to show that for $m = n \cdot 2^{\mu}$ with $\mu \le \nu$
we do not return to the starting number $1$.
For $m = n\cdot 2^0=n$, we have $2^n \equiv 2 \pmod{n}$,
so that $2^n = kn + 2$. For the number $m=n\cdot 2^{\mu}$ we then obtain
\[
2^m=(kn+2)^{2^{\mu}}=(kn)^{2^{\mu}}+2^{\mu}\cdot 2\cdot (kn)^{2^{\mu}-1}+\cdots+
2^{\mu}\cdot 2^{2^{\mu-1}}\cdot (kn)+2^{2^{\mu}}
\]
and therefore $2^m \equiv 2^{2^{\mu}} \pmod{n}$.
But for $0 \le \mu \le \nu$ we have 
$2 \le 2^{2^{\mu}} \le n-1$, and thus 
$2^m \not \equiv 1 \pmod{n}$. For $m = n \cdot 2^{\mu}$ we have therefore
$2^m \not \equiv 1 \pmod{n^2}$.

We have thus shown that for all proper divisors of $n\cdot 2^{\nu+1}$, 
which are smaller than $n\cdot 2^{\nu+1}$, we do not return to the starting number.
Therefore, the invariant set $H_1$ consists of $n\cdot 2^{\nu+1}$ numbers.

In case b) $m=n-1$ is the smallest number for which $g^{m}\equiv 1\pmod{n^2}$ holds.
We choose here as $W_1$ the set of $n-1$ numbers\
$g^{m-1}$, $1\le m\le n-1$, which are calculated modulo $n^2$. 
We can calculate the values
$g^{m-1}$ for any number $m$, $m\ge 1$, but  when
$m= n$, we return to the starting number $1$ and go through the previously obtained values in the same order.

The sets $W_j$ for $j > 1$ we form as follows. 
First, we calculate the\linebreak
 $j$-th number $f_j = \mathrm{rest}(2^{j-1};n^2)$ 
in the invariant set $H_1$ and then we use this number as the starting number for $W_j$. 
Thus, $W_j$ consists of the distinct values $f_jg^{m-1}$ for $m \ge 1$.
Each set $W_j$ consists of $n-1$ distinct numbers from $W$.
Indeed, $f_j \neq 0 \pmod{n}$, $g^{m-1} \neq 0 \pmod{n}$, and therefore also
$f_j g^{m-1} \neq 0 \pmod{n}$. For $m \ge 1$, we have then
\[
f_j g^{m-1} \equiv f_j \Leftrightarrow f_j \cdot (g^{m-1} - 1) \equiv 0
\Leftrightarrow g^{m-1}-1\equiv 0 \Leftrightarrow g^{m-1}\equiv 1\pmod{n^2}
\]
and this means that $W_j$ also consists of $n-1$ distinct elements.
We order the $n-1$ elements in $W_j$ using $f_jg^{m-1}$,
 $m = 1, 2, \cdots, n-1$.
If we allow any number $m$, $m \ge 1$, each set $W_j$ will be cyclically repeated in the same order.

The sets $W_j$ cannot overlap. 
If two sets have a common number, the following values are identical,
and since these sets return to the starting numbers, the sets must be identical. 

For the starting numbers $f_j$ and $f_{j+s}$ of $W_j$ and $W_{j+s}$ with $s>0$, 
we have $f_{j+s} \equiv 2^sf_j \pmod{n^2}$ and $f_j \not \equiv 0 \pmod{n}$. 
Therefore,
\[
W_{j+s}=W_j \Leftrightarrow 2^s f_{j}\in W_j \Leftrightarrow 2^sf_{j}\equiv f_jg^{m}
\Leftrightarrow 2^s\equiv g^m\pmod{n^2}
\]
In the condition $2^s\equiv g^m\pmod{n^2}$, the number $j$ of the first set plays no
role, and this means that all sets must repeat with the same period $s$.
If $j_0$ is the biggest number for which $W_1,\cdots ,W_{j_0}$
are distinct, then in $W_j$ with $1\le j\le n\cdot 2^{\nu+1}$ these sets repeat 
in the same order over and over again. Then
$W_1,\cdots,W_{j_0}$ must obviously contain all numbers from $H_1$
and, due to the cyclic property for the invariant set $H_1$,
$j_0$ must be a divisor of $n\cdot 2^{\nu+1}$.

Now we can see that $j_0=n$. 
It cannot be $j_0=2^{\mu}$ with $2^{\mu}\le 2^{\nu+1}$, because then
\[
j_0\cdot (n-1)=2^{\mu}\cdot (n-1)<n\cdot 2^{\nu+1}
\]
and $W_1,\cdots,W_{j_0}$ do not cover all numbers from $H_1$.
It can't be $j_0=n\cdot 2^{\mu}$ with $\mu\ge 1$ either, because then
\[
j_0\cdot (n-1)=n\cdot 2^{\mu}\cdot (n-1)>n\cdot (n-1)
\]
and $W_1,\cdots,W_{j_0}$ together have more numbers than there are in $W$
and they must overlap. 
Therefore, 
only $j_0=n$ can be true.
The $n$ disjunct sets $W_1,\cdots,W_n$
cover a total of $n(n-1)$ numbers. But as $n(n-1)$ is the quantity of all numbers 
in $W$, the sets $W_j$, $1 \le j \le n$, form a disjunct partition of $W$.

The ordered set $W_j$, $1 \le j \le n$, consists here of the numbers\linebreak
$s_{j,m} = \mathrm{rest}(f_jg^{m-1};n^2)$, $1 \le m \le n-1$, and for $W_j$
we can show directly, analogously to the case of a prime number $n$, that the
numbers in the first and second half of $W_j$ are degrees of inverse elements.
The corresponding set of elements with degrees from $W_j$ consists of inverse pairs, and $x_j$, $1 \le j \le n$, 
\begin{equation}\label {x_j_nn}
 x_j= w^{s_{j,1}}+\cdots+w^{s_{j,n-1}}
\end{equation}
are real numbers.

The polynomial is here $F(x)=(x-x_1)(x-x_2)\cdots(x-x_n)$, and 
for the coefficient $c_p$ of $F(x)$ we obtain the following representation:
\begin{eqnarray}
c_p&=&\mu_{p,0}+\sum_{1\le m\le n-1}\mu_{p,m}w^{n\cdot r_m}
+\sum_{1\le j\le n}\sum_{1\le m\le n-1}\delta_{p,j,m}w^{s_{j,m}}\nonumber 
\end {eqnarray}
In the inner part of the double sum the summands represent the 
elements of the part $W_j$, 
and in the outer part the sum is taken over all possible $W_j$.

If we multiply (we calculate modulo $n^2$) all the numbers in $W$
by the factor $g$,
we obtain a cyclic shift in the summands of $U$ and in the summands of each $W_j$.
This leads to the representation
\begin{eqnarray}
c_p&=&\mu_{p,0}+\sum_{1\le m\le n-1}\mu_{p,m}w^{n\cdot r_{m+1}}
+\sum_{1\le j\le n}\sum_{1\le m\le n-1}\delta_{p,j,m}w^{s_{j,m+1}}
\nonumber 
\end{eqnarray}
where $r_{n}$ is equal to $r_{1}$ and $s_{j,n}$ is equal to $s_{j,1}$.
Thus, we see that here $\mu_{p,m}$,\linebreak
$1\le m\le n-1$, are all equal, and also for every
number $p$ and $j$, the coefficients $\delta_{p,j,m}$, $1\le m\le n-1$, are all equal.
Therefore, the following representation holds for $c_p$:
\begin{eqnarray}
c_p&=&\mu_{p,0}+\mu_{p,1}\cdot\sum_{1\le m\le n-1}v^{r_{m}}
+\sum_{1\le j\le n}\delta_{p,j,1} \sum_{1\le m\le n-1}w^{s_{j,m}}\nonumber \\
&=&\mu_{p,0}-\mu_{p,1}+\sum_{1\le j\le n}\delta_{p,j,1}x_j
\nonumber 
\end{eqnarray}

Next, we multiply (modulo $n^2$) all the numbers in $W$ by the factor $2$.
As a result of this multiplication, $x_1, x_2, \cdots, x_n$ become $x_2, \cdots, x_n, x_1$,
and therefore we obtain the following representation for $c_p$
\begin{eqnarray}
c_p&=&\mu_{p,0}-\mu_{p,1}
+\sum_{1\le j\le n}\delta_{p,j,1} x_{j+1}
\nonumber 
\end{eqnarray}
where $x_{n+1}$ is to be understood as $x_1$. From this we immediately see that
all values of $\delta_{p,j,1}$, $1\le j\le n$, must be equal, and hence
\begin{eqnarray*}
c_p=\mu_{p,0}-\mu_{p,1}
+\delta_{p,1,1}\cdot \sum_{1\le j\le n}x_j=\mu_{p,0}-\mu_{p,1}
+\delta_{p,1,1}\cdot \sum_{m\in W}w^m=\mu_{p,0}-\mu_{p,1}
\end{eqnarray*}
because the sum of the new roots is equal to $0$. 
Thus, $F(x)$ is a polynomial with constant coefficients.

After that we can proceed in case b) exactly as in case a). This completes the Proof 1. 
$\blacksquare$
\section{Square of a Fermat Prime. Proof 2.}
In the second proof of Proposition 2 we do not use the property that the
adapted cyclotomic polynomial is irreducible.
The author believes that it is appropriate to have a proof that does not
use other sophisticated results. 
For this proof we will use the following (somewhat complicated in its formulation)  lemma.
\newline
\newline
{\bf Lemma.} Let $F(x)$ be a polynomial with integer coefficients:
\begin{equation}\label{LemmaPolynom}
F(x)=x^k+c_1x^{k-1}+c_2x^{k-2}+\cdots-c_{k-1}x^1+c_k
\end{equation}
and let the following two properties hold for the roots $x_1,\cdots, x_k$ of $F(x)$:\newline
1) The product $x_ix_j$ of any two roots is a linear combination
of the roots $x_1,\cdots, x_k$ with integer coefficients
and an integer constant:
\begin {eqnarray}
x_ix_j=c_{i,j,1}x_1+c_{i,j,2}x_2+\cdots+c_{i,j,k}x_k +c_{i,j,0}
\end{eqnarray}
2) There exists a prime number $\ell$ such that, for every root $x_j$, 
the coefficients of the representation of $x_j^{\ell}$ as a linear combination
\begin {eqnarray}
x_j^{\ell}=b_{j,1}x_1+b_{j,2}x_2+\cdots +b_{j,k}x_k+b_{j,0}
\end{eqnarray}
satisfy $b_{j,j+1}\equiv 1 \pmod{\ell}$ and $b_{j,m}\equiv 0 \pmod{\ell}$ for 
$m\not= j+1$. 
When $j = k$, the coefficient $b_{k,k+1}$ is taken to be $b_{k,1}$.

Then the following holds true: If one of the roots of $F(x)$ is an integer, then all roots of $F(x)$ are integers.
\newline 
\newline
{\bf Proof.} The fact that $x_j^{\ell}$ can be expressed as an integer constant 
and a linear combination of $x_1,\cdots, x_k$ 
with integer coefficients follows simply from property 1).
We may assume that $x_1=C_1$ is the integer root:
\[
a_{1,1}x_1+a_{1,2}x_2+a_{1,3}x_3+\cdots +a_{1,k}x_k=C_1
\]
where $a_{1,1}=1$, $a_{1,j}=0$ for $j>1$ and $C_1$ is an integer.

Using $x_1^{\ell}$, we obtain the equality
\begin{equation}\label{Gl_2}
a_{2,1}x_1+a_{2,2}x_2+a_{2,3}x_3+\cdots+a_{2,k}x_k=C_2.
\end{equation}
Here, $a_{2,2}\equiv 1\pmod{\ell}$ and $a_{2,j}\equiv 0\pmod{\ell}$ for $j\not=2$.
The value on the right-hand side is $C_2=C_1^{\ell}-b_{1,0}$, and this is an integer.

If we take both sides in
\eqref{Gl_2}
to the power of $\ell$, we obtain the summands $(a_{2,j}x_j)^{\ell}$.
All of them, except $(a_{2,2}x_2)^{\ell}$, are summands with coefficients 
equal to $0$ modulo $\ell$.
With the help of $(a_{2,2}x_2)^{\ell}$ we obtain only one summand which, 
for $x_3$, has a coefficient that is equal to $1$ modulo $\ell$.
For the products obtained by multiplying different terms 
from the left-hand side of
\eqref{Gl_2}
we have coefficients that are equal to $0$ modulo \nolinebreak$\ell$.
This holds for the $\ell$-th power of any sum,
and for the sum on the right-hand side of
\eqref{Gl_2}
it also holds additionally because at least one factor in the corresponding
product has a coefficient that is equal to $0$ modulo $\ell$. 
This leads to the equation
\[
a_{3,1}x_1+a_{3,2}x_2+a_{3,3}x_3+\cdots+a_{3,k}x_k=C_3
\]
where $a_{3,3}\equiv 1\pmod{\ell}$ and $a_{3,j}\equiv 0\pmod{\ell}$ for $j\not=3$.
The value $C_3$ is the difference between $C_2^{\ell}$ and the constants
obtained when calculating the left-hand side. This is again an integer.

We can continue in this way until we get the $k$-th equality. This yields
a system of linear equations with integer coefficients and integer values
on the right-hand side
\begin{eqnarray*}
&&x_1\;\;\;\; +\;\;0x_2\;\;\;+\;\;0x_3\;+\cdots +\;\;0x_k\;=C_1\\
&&a_{2,1}x_1+a_{2,2}x_2+a_{2,3}x_3+\cdots+a_{2,k}x_k=C_2\\
&&a_{3,1}x_1+a_{3,2}x_2+a_{3,3}x_3+\cdots+a_{3,k}x_k=C_3\\
&&\cdots\\
&&a_{k,1}x_1+a_{k,2}x_2+a_{k,3}x_3+\cdots+a_{k,k}x_k=C_k
\end{eqnarray*}
For the coefficients, we have $a_{i,i}\equiv 1\pmod{\ell}$
and $a_{i,j}\equiv 0\pmod{\ell}$ for $i\not=j$.

The first term $a_{j,1}x_1$ of the $j$-th equation, $2 \le j \le k$,
can be eliminated by subtracting the first equation,
multiplied by the factor $a_{j,1}$, from the $j$-th equation. This yields the system of equations
\begin{eqnarray*}
&&x_1\;+\;0x_2\;\;\;+\;\;0x_3+\cdots +\;\;0x_k\;=C_1\\
&&0x_1+a'_{2,2}x_2+a'_{2,3}x_3+\cdots+a'_{2,k}x_k=C'_2\\
&&0x_1+a'_{3,2}x_2+a'_{3,3}x_3+\cdots+a'_{3,k}x_k=C'_3\\
&&\cdots\\
&&0x_1+a'_{k,2}x_2+a'_{k,3}x_3+\cdots+a'_{k,k}x_k=C'_k
\end{eqnarray*}
The coefficients of this system and the values on the right-hand side
are integers.
Here, 
the coefficients satisfy 
$a'_{i,i}\equiv 1\pmod{\ell}$
and $a'_{i,j}\equiv 0\pmod{\ell}$ when $i \neq j$.

For $3 \le j \le k$, we can first multiply the $j$-th equation by the factor 
$a'_{2,2}$ ($a'_{2,2} \neq 0$) and then subtract the second equation, 
multiplied by the factor $a'_{j,2}$, from the resulting equation.
This leads to the system of equations
\begin{eqnarray*}
&&x_1\;+\;0x_2\;\;\;+\;\;0x_3+\cdots +\;\;0x_k\;=C_1\\
&&0x_1+a'_{2,2}x_2+a'_{2,3}x_3+\cdots+a'_{2,k}x_k=C'_2\\
&&0x_1+0x_2\;\;\;+a''_{3,3}x_3+\cdots+a''_{3,k}x_k=C''_3\\
&&\cdots\\
&&0x_1+0x_2\;\;\;+a''_{k,3}x_3+\cdots+a''_{k,k}x_k=C''_k
\end{eqnarray*}
The coefficients of this system and the values on the right-hand side 
are integers. The coefficients on the main diagonal are again
equal to $1$ modulo \nolinebreak$\ell$, and the others are equal 
to $0$ modulo $\ell$.

In this way, we  ultimately get a system of linear equations
with integer coefficients and integer values on the right-hand side,
in which all coefficients below the main diagonal are equal to $0$
and the coefficients on the main diagonal are not equal to $0$.
The solutions of this system are obviously rational numbers. 
However, since a rational root of the polynomial
\eqref{LemmaPolynom}
must necessarily be an integer,   
all roots $x_1,\cdots,x_k$ are integers. This completes the proof.
$\blacksquare$
\newline
\newline
{\bf Proof 2.} In this proof of Proposition 2 we continue to use 
the polynomial $F(x)=(x-x_1)(x-x_2)\cdots(x-x_n)$ of odd degree $n$, 
for which we already know that it is a polynomial with integer coefficients.
We can assume here that $g$ is a prime number. Indeed, we have 
$g^{\frac{n-1}{2}}\not\equiv 1 \pmod{n}$, and $\frac{n-1}{2}$ is the greatest 
divisor of $n-1=2^{2^{\nu}}$. If $g_j$, $1\le j\le j_0$, are the prime divisors of $g$,
then for at least one of them, $g_j^{\frac{n-1}{2}}\not\equiv 1\pmod{n}$
must hold, and $g_j$ is the appropriate primitive root.
We set $\ell=g$ in case a) and $\ell=2$ in case \nolinebreak b) and have to
show that the  conditions of the lemma proven above hold for the roots 
$x_1,\cdots,x_n$ of $F(x)$.

If we calculate the product $x_{j_1}x_{j_2}$ of two roots, we can obtain an integer constant
and elements with degrees from $U$ and from $W_j$, $1\le j\le n$,
\[
x_{j_1}x_{j_2}=\xi_{j_1,j_2,0}+\sum_{1\le m\le n-1}\xi_{j_1,j_2,m}w^ {n\cdot r_m}
+\sum_{1\le j\le n}\sum_{1\le m\le n-1}\theta_{j_1,j_2,j,m}w^{s_{j,m}}\nonumber 
\]
If we multiply (modulo $n^2$ calculated) the numbers in $W_{j_1}$ and $W_{j_2}$ 
by the factor $g^n$ in case a) and by the factor $g$ in case b), 
then $x_{j_1}$ and $x_{j_2}$ remain unchanged, and we obtain the representation
\[
x_{j_1}x_{j_2}=\xi_{j_1,j_2,0}+\sum_{1\le m\le n-1}\xi_{j_1,j_2,m}w^{n\cdot r_{m+1}}
+\sum_{1\le j\le n}\sum_{1\le m\le n-1}\theta_{j_1,j_2,j,m}w^{s_{j,m+1}}\nonumber 
\]
Thus, we see that for the elements with degrees from $U$ and for 
the elements with degrees from the parts $W_j$ we have even coverings
and obtain  the required representation
\begin{eqnarray*}
x_{j_1}x_{j_2}&=&\xi_{j_1,j_2,0}+\xi_{j_1,j_2,1}\sum_{1\le m\le n-1}v^{r_m}
+\sum_{1\le j\le n}\theta_{j_1,j_2,j,1}\sum_{1\le m\le n-1}w^{s_{j,m}}\nonumber \\
&=&\xi_{j_1,j_2,0}-\xi_{j_1,j_2,1} +\sum_{1\le j\le n}\theta_{j_1,j_2,j,1}x_j
\end{eqnarray*}

If we calculate $x_j^{\ell}$, we obtain, because of 
\eqref{x_j_nn},
the sum of the elements of degree $\ell$ for the individual summands of $x_j$,
and that is the value $x_{j+1}$
\[
\sum_{m\in W_j}w^{\ell\cdot m}=\sum_{m\in W_{j+1}}w^m=x_{j+1},
\]
where $x_{n+1}$ is to be understood as $x_1$.
In addition to $x_{j+1}$ we obtain here for the summands $w^{s_{j,1}},\cdots,w^{s_{j,n-1}}$
of $x_j$ the corresponding products with coefficients divisible by $\ell$.
The elements obtained in this way, and therefore also the values $x_j$ 
calculated in this way, all have coefficients that are divisible by $\ell$. 
The second condition of the lemma is then satisfied.
\newline

Now we want to consider the following values $d_p$, $1 \le p \le n$,
\begin{eqnarray}\label{dp}
d_p = x_1^p + x_2^p + \cdots + x_n^p
\end {eqnarray}
These values are easier to analyze than $c_p$, and we will
see later that they will be helpful for analyzing the coefficients $c_p$. 
Here we have $d_1=0$, and we consider the values $d_p$, $2\le p\le n$.

Let us first assume $2 \le p < n$. 
If we consider just one summand $x_j^p$, we see from
\eqref{x_j_nn}
that we obtain terms $w^{m p}$ for $m\in W_j$ and terms with the inner products
for $w^{m_1},w^{m_2},\cdots$.
Since $m\in W_j$, it follows that $m p \in W$, and we have the following representation 
for $x_j^p$:
\begin{eqnarray}\label {x_j^p_nn}
x_j^p=\theta_{p,0}
+\{\tau_{p,1}w^{1\cdot n}+\cdots+ \tau_{p,n-1}w^{(n-1)\cdot n}\}+\cdots
\end{eqnarray}
The summands with degrees from $W$ are not shown in
\eqref{x_j^p_nn}
in more detail. For $x_j^p$, we have the constant $\theta_{p,0}$  and
$\theta_{p,1}=\tau_{p,1}+\cdots+ \tau_{p,n-1}$ 
elements with degrees from $U$, where $\theta_{p,0}$ and $\theta_{p,1}$ are integers. 

In case a) we multiply all the numbers in $W$ by the factor $g$.
In this case $x_1, x_2, \cdots, x_n$ become $x_2, \cdots, x_n, x_1$,
and instead of 
\eqref{x_j^p_nn}
for $x_j^p$, we obtain the representation for $x_{j+1}^p$ as
\begin{eqnarray}\label{x_j+1^p_nn}
x_{j+1}^p=\theta_{p,0}
+\{\tau_{p,1}w^{ n\cdot g}+\cdots+ \tau_{p,n-1}w^{(n-1)\cdot n\cdot g}\}+\cdots
\end{eqnarray}
and we see that in this case we also have $\theta_{p,1}$ elements with degrees from $U$.
If we then sum over all $j$, $1 \le j \le n$, we obtain a total of
$n\cdot \theta_{p,1}$ elements with degrees from $U$. The height of the even coverage 
of $U$ is therefore equal to $n\cdot \theta_{p,1}/(n-1)$,
and for $d_p$, $2\le p<n$, the following holds: 
\begin{eqnarray}\label {dp_in_U}
d_p=n\cdot \theta_{p,0}-n\cdot \theta_{p,1}/(n-1)
=n\cdot\mu_{p,0},
\end{eqnarray}
where $\mu_{p,0}$ is an integer.

In case b), for $2 \le p < n$, we go from $x_1,x_2,\cdots, x_n$ to
$x_2, \cdots, x_n, x_1$ by multiplying all numbers in $W$ by the factor of $2$.
Using this multiplication we obtain for $d_p$ in case b) the same representation
\eqref{dp_in_U}.

The procedure for $p=n$ is analogous. But, instead of
\eqref{x_j^p_nn} and \eqref{x_j+1^p_nn},
we obtain in this case the presentations
\begin{eqnarray}\label {x_j^n_nn}
x_j^n=n\cdot \theta_{n,0}
+n\cdot\{\tau_{n,1}w^{1\cdot n}+\cdots+ \tau_{n,n-1}w^{(n-1)\cdot n}\}+\cdots
\end{eqnarray}
and
\begin{eqnarray}\label {x_{j+1}^n_nn}
x_{j+1}^n=n\cdot\theta_{n,0}
+n\cdot\{\tau_{n,1}w^{ n\cdot g}+\cdots+ \tau_{n,n-1}w^{(n-1)\cdot n\cdot g}\}+\cdots
\end{eqnarray}
since the coefficients of the corresponding inner products must all be divisible by the prime number $n$.
For any $m\in W_j$ in this case we have $n\cdot m\in U$, 
and for each of the $n-1$ numbers in $W_j$ we obtain a number from $U$.
Calculating the sum over all $j$, $1\le j\le n$, we obtain then a total of
$n(n-1)$ elements with degrees from $U$. They therefore provide  
an additional even coverage of the height $n=n(n-1)/(n-1)$
for the elements with degrees from $U$.
Hence, for $p=n$, instead of
\eqref{dp_in_U},
we obtain the following representation:
\begin{eqnarray}\label{dn_in_U}
d_n=n^2\cdot\theta_{n,0}-n^2 \theta_{n,1}/(n-1) -n =n^2\cdot\mu_{n,0} -n
\end{eqnarray}

The coefficient $c_1$ of the polynomial $F(x)$ is known here, $c_1=0$, 
and for the coefficient $c_p$ the following representation holds: 
$c_p = (-1)^p\hat{c}_p$, 
where
\[
\hat{c}_p=\sum_{M(n,p)}x_{j_1}x_{j_2}\cdots x_{j_p}
\]
and $\hat{c}_1=-c_1$ is here equal to $0$.

We are going to examine the relationship between the values $\hat{c}_p$ and $d_p$.
To do this, we need a simple generalization for the sets $M(n,p)$ used above.  
For $n$, $p$ and $s$, $2\le s\le p\le n$, we denote
$M(n,p,s)$ the set of the sequences $\{j_1,j_2,\cdots,j_p\}$ of integers,
for which $1\le j_2 \le \cdots\le j_p\le n$ and which are formed so that
exactly $s$ of the numbers $j_1,j_2,\cdots, j_p$ are equal.

For $p \ge 2$ we have 
\begin{eqnarray}\label{cp_kalk_1}
\hat{c}_{p-1}\cdot d_1 =
\left\{ \sum_{M(n,p-1)} x_{j_1}x_{j_2}\cdots x_{j_{p-1}} \right\}
\cdot (x_1+x_2+\cdots +x_n)\;\;\;\;\;\;\;\;\;\;\;\;\;\;\;\;\;\;\;\;\;\;\;\;\;\;\; \nonumber \\
=p\sum_{M(n,p)}x_{j_1}x_{j_2}\cdots x_{j_p}
-\sum_{M(n,p,2)}x_{j_1}x_{j_2}\cdots x_ {j_p}
=p\cdot \hat{c}_p-\sum_{M(n,p,2)}x_{j_1}x_{j_2}\cdots x_{j_p}\;\;
\end{eqnarray}
and the last sum in
\eqref{cp_kalk_1}
can be represented as follows:
\begin{eqnarray*}
&&\sum_{M(n,p,2)}x_{j_1}x_{j_2}\cdots x_{j_p}
=\sum_{1\le j\le n}x_j^2
\sum_{M(n,p-2)}x_{j_1}x_{j_2}\cdots x_{j_{p-2}}
-\sum_{M(n,p,3)}x_{j_1}x_{j_2}\cdots x_{j_p}\;\;\;\;\;\;\;\;\;\;\;\;\;\;\;\;\;\;\;\;\;\;
\;\;\;\;\; \\
&&=\hat{c}_{p-2}\cdot d_2 -\sum_{M(n,p,3)}x_{j_1}x_{j_2}\cdots x_{j_p}
\end{eqnarray*}
Similarly, for $3 \le s \le p-1$ we obviously obtain the equality
\[
\sum_{M(n,p,s)}x_{j_1}x_{j_2}\cdots x_{j_p}=\hat{c}_{p-s}\cdot d_s -\sum_{M(n,p,s+1)}x_{j_1}x_{j_2}\cdots x_{j_p}
\]
Taking these results into account, and noting that  
\[
\sum_{M(n,p,p)}x_{j_1}x_{j_2}\cdots x_{j_p}=\sum_{1\le j\le n} x_j^p=d_p,
\]
we obtain, based on 
\eqref{cp_kalk_1},
the representation:
\begin{equation}\label{|cp|_dp}
p\cdot \hat{c}_p= \hat{c}_{p-1}d_1-\hat{c}_{p-2}d_2+\hat{c}_{p-3}d_3+\cdots+(-1)^{p-1} d_p.
\end{equation}
The term $\hat{c}_{p-1}d_1$ in this formula is equal to $0$, and is only
shown to make it look simpler. 

We have here $\hat{c}_1=0$, and for the values $d_p$ it is clear that 
$d_p\equiv 0\pmod{n}$ for $1\le p<n$ and $d_n/n\equiv-1\pmod{n}$. 
Using these properties, based on 
\eqref{|cp|_dp},
we obtain for $\hat{c}_p$, when $1\le p<n$,
\[
\hat{c}_p\equiv\frac{1}{p}( \hat{c}_{p-1}\cdot 0-\hat{c}_{p-2}\cdot 0+\hat{c}_{p-3}\cdot 0+\cdots+(-1)^{p-1}\cdot 0)\equiv 0\pmod{n}.
\]
For $\hat{c}_n$ we obtain
\[
\hat{c}_n\equiv \frac{1}{n}(\hat{c}_{n-1}\cdot 0-\hat{c}_{n-2}\cdot 0+\hat{c}_{n-3}\cdot 0+\cdots -\hat{c}_1\cdot 0)
-\frac{d_n}{n} \equiv -1\pmod{n}.
\]
For the coefficients $c_p$ we have therefore the following: $c_p\equiv 0\pmod{n}$ for $1\le p<n$
and $c_n\equiv 1\pmod{n}$.

If $x$ is an integer solution of the polynomial $F(x)$, then  
\[
x^n+c_2 x^{n-2}+\cdots+c_{n-1}x+c_n\equiv x^n+0\cdot x^{n-2}+\cdots +0\cdot x +1
\equiv x^n+1\equiv 0\pmod{n}
\]
Due to $x\not \equiv 0\pmod{n}$ we have here $x^n\equiv x\pmod{n}$ 
and therefore\linebreak
 $x\equiv -1\pmod{n}$.
For every root $x_j$ of $F(x)$ we have thus the representation 
$x_j=-1+k_j\cdot n$. 
But $x_j$ is the sum of $(n-1)/2$ double cosine values, and therefore $-n < x_j < n$. All roots 
$x_1,\cdots,x_n$ must thus be equal to $-1$.
This leads to a contradiction, because then $x_1+x_2+\cdots +x_n=-n$ holds, and not 
$x_1+x_2+\cdots +x_n=0$. 
Since $F(x)$ has no integer solution, it is not possible to construct all values of $x_1,\cdots,x_n$.
Therefore the corresponding regular polygon cannot be constructed using a compass and straightedge.
This concludes the Proof 2.
$\blacksquare$

{\bf Remark 1.} Unfortunately, the proof given here cannot simply
be extended to the case of primes that are not Fermat numbers. 
The lemma is applicable, and it is also possible to represent $\hat{c}_m$
in terms of $d_m$. However, for the values $d_m$, and therefore
also for $\hat{c}_m$, there is no such simple calculation that allows the transition 
to an equation for the values of $\mathrm{rest}(x_j; k)$.

{\bf Remark 2.} The proof for squares of Fermat primes 
is somewhat more complicated 
than for primes that are not 
Fermat numbers. Two cases must be considered here, and in case b) it must be 
understood how to form the appropriate disjunct partitions of $W$.
For all currently known Fermat primes, we have the case a).
However, since the existence of other Fermat primes and the corresponding
property of their squares is unknown, the case b) could not be omitted from the proof. 

\section{Examples}

The proof of Wantzel’s theorem presented here has the advantage
that it shows the way in which we can for specific (yet interesting) cases
simply demonstrate that the corresponding regular $n$-gon
cannot be constructed using a compass and straightedge. Two things are important here:

1) The property that the proposed approach yields a polynomial
with integer coefficients. 

2) The property that not all real roots of a polynomial of odd degree 
with integer coefficients can be constructed if it has no integer roots.

For values of $n$ that are not too large, the approximations for the coefficients 
of the corresponding polynomial can be calculated numerically.
To do this, one can, for example, successively calculate the coefficients 
of the polynomials $x-x_1$, $(x-x_1)(x-x_2)$,...,
$(x-x_1)(x-x_2)\cdots (x-x_k)$ numerically.
Using the approximations for the coefficients, the corresponding
polynomial with integer coefficients can be determined.
For the polynomial obtained in this way it is easy to show 
that it has no integer roots. The more complicated steps from the proof 
outlined above are not necessary in specific cases.

We know that every value $x_j$ is the sum of pairs of inverse elements,
and when we order their numbers, we can proceed in such a way that we get the numbers of the corresponding pairs. 
There are half as many numbers of pairs as numbers of elements.
The number of the next element, based on the number $s$ of the previous one, is
$\mathrm{rest}(s\cdot g;n)$.
If we want to 
determine the number of the next pair we have to calculate
$\min(\mathrm{rest}(s\cdot g;n);n-\mathrm{rest}(s\cdot g;n))$.
Here, if desired, we can analytically calculate the products of pairs using
the equality
$p_{k_1}\cdot p_{k_2}=p_{|k_1-k_2|}+p_{\min(k_1+k_2,n-(k_1+k_2))}$,
\cite{Mainik}.
\newline
\newline
{\bf Example 1: Regular $7$-gon.} For $n=7$, $g=3$ is a primitive root.
The ordered numbers of pairs are here $1,3,2$. 
From $2$ we return to the starting number $1$.
Here we have $n-1=6=2\cdot 3$, and the corresponding $3$ values:
$x_1=p_1$, $x_2=p_3$, and $x_3=p_2$.
In this simple example the resulting polynomial can be determined analytically.
The coefficients of this polynomial are:\newline
$c_1=-(p_1+p_3+p_2)=1$\newline
$c_2=p_1p_3+p_1p_2+p_3p_2= p_2+p_3+p_1+p_3+p_1+p_2=2(p_1+p_2+p_3)=-2$\newline
$c_3=-p_1p_3p_2=-(p_2+p_3)p_2=-(p_3+2+p_1+p_2)=-(2-1)=-1$\newline
and the polynomial is $x^3+x^2-2x-1$. This polynomial has no integer roots.
Only the divisors of $c_3=-1$ are possible, and $1$ and $-1$ are not roots.
The roots $x_1$, $x_2$, $x_3$ cannot be constructed with a compass and straightedge, and
therefore the regular $7$-gon cannot be constructed as well.

The polynomial can also be determined numerically. We have here:\linebreak
$x_1=2\cos(1\cdot 2\pi/7)\approx 1.246979604$,
$x_2=2\cos(3\cdot 2\pi/7)\approx -1.801937736$,
$x_3=2\cos(2\cdot 4\pi/7) \approx -0.445041868$,
and the coefficients are obtained numerically:
$c_1\approx 1.000000$, $c_2\approx -2.000000$, $c_3\approx -1.000000$.
These approximations for the coefficients show that the corresponding
polynomial (with integer coefficients) is $x^3+x^2-2x-1$.
This is the same polynomial that we previously calculated analytically,
and it doesn't have integer roots.

{\bf Remark 1.} We have also calculated the coefficients $c_1$ and $c_2$ here,
in order to show the procedure in more detail.

{\bf Remark 2.} We could already understand from the approximate numerical values $x_j$
that we will get a polynomial with roots that are close to the values $x_j$
and therefore cannot be integers. However, this two-step estimation 
(first for the coefficients and then for the new roots)
is uninteresting and unnecessary.
It is enough to determine the correct polynomial with integer coefficients and directly
establish that this polynomial has no integer roots.
\newline
\newline
{\bf Example 2: Regular $13$-gon.} For $n=13$ is $g=2$ a primitive root.
The ordered numbers of the pairs are $1, 2, 4, 5, 3, 6$.
Here we have $n-1=4\cdot 3$ and calculate 
$x_1,x_2,x_3$ as follows:\newline
$x_1=2\cos(1\cdot 2\pi/13)+2\cos(5\cdot 2\pi/13)\approx 0.2738905549$,\newline
$x_2=2\cos(2\cdot 2\pi/13)+2\cos(3\cdot 2\pi/13)\approx 1.377202853$,\newline
$x_3=2\cos(4\cdot 4\pi/13)+2\cos(6\cdot 2\pi/13)\approx -2.651093409$. \newline
For the coefficients of the corresponding polynomial we obtain $c_1\approx 0.999999$,\linebreak
$c_2 \approx -4.000000$, $c_3 \approx 1.000000$,
and the polynomial is
$F(x) = x^3 + x^2 - 4x + 1$.
It is also immediately apparent that the polynomial
$F(x)$ does not have an integer root.
The roots $x_1, x_2, x_3$, and therefore also the regular $13$-gon, cannot be constructed
with a compass and straightedge.
\newline
\newline
{\bf Example 3: Regular $9$-gon.} $9$ is the square of the Fermat prime $3$.
For $n=3$, the number $g=2$ is a primitive root.
The ordered numbers of the pairs for the new roots are here:
$1, 2, 4$. 
We have only $3$ values here: $x_1=p_1$, $x_2=p_2$, $x_3=p_4$.
In this simple example the coefficients can be calculated analytically: 
\newline
\newline
$c_1=-(p_1+p_2+p_3)=0$,\newline
$c_2=p_1p_2+p_1p_4+p_2p_4=p_1+p_3+p_3+p_4+p_2+p_3=
(p_1+p_2+p_4)+3p_3=-3$\newline
$c_3=-p_1p_2p_4=-(p_1+p_3)p_4=-(p_3+p_4+p_1+p_2)=-p_3=1$,\newline
and the corresponding polynomial is $x^{3}-3x+1$.
This polynomial has no integer roots. Therefore the values $x_1,x_2,x_3$
and thus also the regular $9$-gon cannot be constructed with a compass and straightedge.

Here we can also calculate numerically:
$x_1=2\cos(1\cdot 2\pi/9)\approx 1.532088886$, 
$x_2=2\cos(2\cdot 2\pi/9\approx 0.347296355$,  
$x_3=2\cos(4\cdot 2\pi/9)\approx -1.879385242$.
Based on these approximations, we obtain 
the coefficients $c_1\approx 0.000000$, $c_2\approx -3.000000$,
$c_3\approx 1.000000$ and again the polynomial $x^{3}-3x+1$.

The presented approaches
can be interpreted as slightly different proofs that
trisection with a compass and straightedge is not possible for every angle.
\newline
\newline
{\bf Example 4: Regular $25$-gon.} $25$ is the square of the Fermat prime $5$.
For $n=5$ the number $g=3$ is a primitive root.
Here the ordered numbers of the pairs are $1,3,9,2,6,7,4,12,11,8$. 
Based on these $10$ pairs, we calculate $5$ values $x_j$ adding the next
$5$-th pair to the $j$-th pair:\newline
$x_1=2\cos(1\cdot 2\pi/25)+2\cos(7\cdot 2\pi/25)\approx 1.562403693$,\newline
$x_2=2\cos(3\cdot 2\pi/25)+2\cos(4\cdot 2\pi/25)\approx 2.529590845$,\newline
$x_3=2\cos(9\cdot 4\pi/25)+2\cos(12\cdot 2\pi/25)\approx -3.259077382$,\newline
$x_4=2\cos(2\cdot 2\pi/25)+2\cos(11\cdot 2\pi/25)\approx -0,106939612$, \newline
$x_5=2\cos(6\cdot 4\pi/25)+2\cos(8\cdot 2\pi/25)\approx -0,725977544$.\newline
By numerically calculation we obtain the following approximations for the coefficients
of the corresponding polynomial:
$c_1\approx 0,000000$, $c_2\approx -9,999999$, $c_3\approx 5,000000$,
$c_4\approx 9,999999$, $c_5\approx 0,999999$.
The polynomial (with integer coefficients) is therefore $x^5-10x^3+5x^2+10x+1$, and
it is easy to see that this polynomial does not have an integer root.
Not all roots $x_1, \cdots, x_5$ can be constructed, and therefore the regular $25$-gon cannot be constructed with a compass and straightedge.
\newline
\newline
{\bf Example 5: Regular $17^2$-gon.} In this case the approximations for the values $x_1,x_2,\cdots, x_{17}$ 
and for the coefficients of the corresponding polynomial still can be easily calculated numerically. Even an Excel spreadsheet would suffice here.
This leads to the following polynomial (with integer coefficients):\newline
$x^{17}-136x^{15}+85x^{14}+6154x^{13}-6545x^{12}-119680x^{11}+168555x^{10}
+998835x^{9}\linebreak
-1749300x^{8}-2783546x^{7}+6581040x^{6}-678725x^{5}-3813882x^{4}
+770593x^{3}\linebreak
+616267x^{2}-82620x-577$\newline
This polynomial does not have an integer root, 
and thus the regular $289$-gon cannot be constructed with a compass and straightedge.

\end{document}